\documentclass{amsart}

\usepackage{amsthm, amsfonts, amssymb, amsmath, latexsym, enumerate, times}
\usepackage[latin1]{inputenc}
\usepackage{mathrsfs}
\usepackage{mathtools}
\usepackage{letltxmacro}
\LetLtxMacro\orgvdots\vdots
\LetLtxMacro\orgddots\dots

\makeatletter
\DeclareRobustCommand\vdots{%
	\mathpalette\@vdots{}%
}
\newcommand*{\@vdots}[2]{%
	\sbox0{$#1\cdotp\cdotp\cdotp\m@th$}%
	\sbox2{$#1.\m@th$}%
	\vbox{%
		\dimen@=\wd0 %
		\advance\dimen@ -3\ht2 %
		\kern.5\dimen@
		\dimen@=\wd2 %
		\advance\dimen@ -\ht2 %
		\dimen2=\wd0 %
		\advance\dimen2 -\dimen@
		\vbox to \dimen2{%
			\offinterlineskip
			\copy2 \vfill\copy2 \vfill\copy2 %
		}%
	}%
}
\DeclareRobustCommand\ddots{%
	\mathinner{%
		\mathpalette\@ddots{}%
		\mkern\thinmuskip
	}%
}
\newcommand*{\@ddots}[2]{%
	\sbox0{$#1\cdotp\cdotp\cdotp\m@th$}%
	\sbox2{$#1.\m@th$}%
	\vbox{%
		\dimen@=\wd0 %
		\advance\dimen@ -3\ht2 %
		\kern.5\dimen@
		\dimen@=\wd2 %
		\advance\dimen@ -\ht2 %
		\dimen2=\wd0 %
		\advance\dimen2 -\dimen@
		\vbox to \dimen2{%
			\offinterlineskip
			\hbox{$#1\mathpunct{.}\m@th$}%
			\vfill
			\hbox{$#1\mathpunct{\kern\wd2}\mathpunct{.}\m@th$}%
			\vfill
			\hbox{$#1\mathpunct{\kern\wd2}\mathpunct{\kern\wd2}\mathpunct{.}\m@th$}%
		}%
	}%
}
\makeatother

\usepackage{array}
\usepackage{comment}
\usepackage{paralist}
\usepackage{xcolor}

\newtheorem{theorem}{Theorem}

\newtheorem{conjecture}[theorem]{Conjecture}

\theoremstyle{definition}

\newtheorem{remark}[theorem]{Remark}

\newcommand{\bbC}{{\mathbb C}}

\def\le{\leqslant}
\def\ge{\geqslant}

\begin{document}
 
\title[A remark on the canonical degree]{A remark on the canonical degree \\ of curves on smooth projective surfaces}
 
\author{Ciro Ciliberto}
\address{Dipartimento di Matematica, Universit\`a di Roma Tor Vergata, Via O. Raimondo 00173 Roma, Italia}
\email{cilibert@axp.mat.uniroma2.it}

\author{Claudio Fontanari}
\address{Dipartimento di Matematica, Universit\`a degli Studi di Trento, Via Sommarive 14, 38123 Povo, Trento}
\email{claudio.fontanari@unitn.it}
 
 

\centerline{}
\begin{abstract} 
The canonical degree $C.K_X$ of an integral curve on a smooth projective surface $X$ is conjecturally bounded 
from above by an expression of the form $A(g-1)+B$, where $g$ is the geometric genus of $C$ and $A$, $B$ 
are constants depending only on $X$. We prove that this conjecture holds with $A = -1$ under the assumptions 
$h^0(X, -K_X) = 0$ and $h^0(X, 2K_X + C) = 0$. 
\end{abstract}

\maketitle 

\section{Introduction} 
Let $C$ be an integral curve on a smooth complex projective surface $X$. We denote by $g = g(C)$ 
its geometric genus and by by $k_C := C.K_X$ its \emph{canonical degree}. 
The following conjecture (see \cite{CR}, Conjecture 5.1) is widely open: 

\begin{conjecture} \label{main} 
Let $X$ be a smooth projective surface. There exist constants $A$, $B$ such that for any integral curve $C$ we have
$$k_C \le A(g - 1) + B.$$
\end{conjecture}

Notice that if $h^0(X, -K_X) > 0$ then $k_C \le 0$ for all but finitely many curves $C$ on $X$, hence we
can assume $h^0(X, -K_X) = 0$. 

Moreover, if $C$ is smooth then $k_C = 2(g-2) - C^2$, so in this case a lower bound for $C^2$ in terms of 
$g$ provides an answer to Conjecture \ref{main}. In particular, Hao's proof of the so-called \emph{Weak 
bounded negativity conjecture} yields the following result:

\begin{theorem}\label{Hao} 
Let $X$ be a smooth projective surface with $h^0(X, -K_X) = 0$ and let $C$ be a smooth irreducible 
curve of genus $g$ on $X$. 

(a) If  $h^0(X, 2K_X + C) \ne 0$  then $k_C \le  4(g - 1) + 3 c_2(X) - K_X^2$. 

(b) If   $h^0(X, 2K_X + C) = 0$  then $k_C \le -(g-1) - K_X^2 - \chi (\mathcal{O}_X)$. 

\end{theorem}

Item (a) is precisely \cite{H}, Corollary 1.8  (just notice that $h^0(X, 2K_X + C) \ne 0$ implies 
$h^0(X, 2(K_X + C)) \ne 0$),  while item (b) follows from the proof of \cite{H}, Lemma 1.3, 
which  works under the weaker assumption $h^0(X, 2K_X + C) = 0$ and  in the smooth case 
gives the better bound $C^2 \ge 3 g + K_X^2 + \chi (\mathcal{O}_X) - 3$. 

The extension of item (a) to singular curves is known to be a highly nontrivial problem. A partial result 
towards this direction is \cite{LM}, Theorem 1, item (4): if $C$ is integral and some multiple of $K_X + C$ 
is effective then $k_C \le  4(g - 1) + 3 c_2(X) - K_X^2 + n$, where $n$ is the number of ordinary nodes 
and ordinary triple points of $C$. As pointed out by Miyaoka in \cite{M}, Remark D., this yields \textit{a similar 
bound for $C.K_X$ provided $C$ contains neither ordinary double points nor ordinary triple points.
Strangely, complicated singularities of high multiplicity do no harm to estimating
canonical degrees. Curves with many ordinary double points are technically
the most difficult to deal with.}

Here we focus on item (b) and we obtain the following complete generalization to the singular case: 

\begin{theorem}\label{singular}
Let $X$ be a smooth projective surface with $h^0(X, -K_X) = 0$ and let $C$ be an integral 
curve of geometric genus $g$ on $X$. If $h^0(X, 2K_X + C) = 0$ then 
$$
k_C \le -(g-1) - K_X^2 - \chi (\mathcal{O}_X).
$$
\end{theorem}

Our argument is inspired by \cite{H}, proof of Theorem 1.9, which in turn closely follows \cite{BB+}, proof of Proposition 3.5.3. 

We work over the complex field $\bbC$.

\medskip

{\bf Acknowledgements:} The authors are members of GNSAGA of the Istituto Nazionale di Alta Matematica ``F. Severi". This research project was partially supported by PRIN 2017 ``Moduli Theory and Birational Classification''. 

\section{The proofs}

\noindent \textit{Proof of Theorem \ref{singular}.} The idea is to blow up $X$ resolving step by step the singularities of $C$. 
If the assumptions hold at each step provided they hold at the previous one and the conclusion holds at each step provided 
it holds at the next one, then the statement follows recursively from the smooth case, namely from item (b) of Theorem \ref{Hao}. 

Let $X$ be a smooth projective surface and let $C$ be an integral curve of geometric genus $g$ on $X$. 
Let $p \in C$ be a point with multiplicity $\mathrm{mult}_p(C) = m \ge 2$ and let $\pi: \tilde{X} \to X$ 
be the blow-up of $X$ at $p$.
Let $E$ be the exceptional divisor of the blow-up and let $\tilde{C} = \pi^*(C)-mE$ be the strict transform of $C$. 

We claim the following: 

(i) If $h^0(X, -K_X) = 0$ then $h^0(\tilde{X}, -K_{\tilde{X}}) = 0$.

 (ii) If $h^0(X, 2K_X + C) = 0$ then $h^0(\tilde{X}, 2K_{\tilde{X}} + \tilde{C}) = 0$.
 
(iii) If $k_{\tilde{C}} \le -(g-1) - K_{\tilde{X}}^2 - \chi (\mathcal{O}_{\tilde{X}})$ then 
$k_C \le -(g-1) - K_X^2 - \chi (\mathcal{O}_X)$.

Indeed, (i) follows from $K_{\tilde{X}} = \pi^*(K_X) + E$ (for details see \cite{H}, Lemma 1.10). 

Next, for (ii) we have $h^0(\tilde{X}, 2K_{\tilde{X}} + \tilde{C}) = h^0(\tilde{X}, 2 (\pi^*(K_X) + E) + \pi^*(C)-mE)) \le 
h^0(\tilde{X}, 2 \pi^*(K_X) + \pi^*(C)) = h^0(\tilde{X}, \pi^*(2K_X + C))  = h^0(X, 2K_X + C) = 0$. 
 
Finally, we have $k_C = K_X . C = K_{\tilde{X}} . \tilde{C} - m = k_{\tilde{C}} - m \le 
-(g-1) - K_{\tilde{X}}^2 - \chi (\mathcal{O}_{\tilde{X}}) - m = 
-(g-1) - (K_X^2 - 1) - \chi (\mathcal{O}_X) - m \le -(g-1) - K_X^2 - \chi (\mathcal{O}_X).$

\qed

\begin{remark}
If the curve $C = \sum_i^ n C_i$ is reducible (but still reduced), the argument above works verbatim by setting 
$g :=  \sum_i^ n g_i - (n-1)$, where $g_i$ is the geometric genus of the irreducible component $C_i$. 
Indeed, after finitely many blow-ups we obtain a curve $\tilde{C}$ which is the disjoint union of the 
normalizations of the curves $C_i$. The arithmetic genus of $\tilde{C}$ is $p_a(\tilde{C}) = 
\sum_i^ n g_i - (n-1) = g$ and the proof of \cite{H}, Lemma 1.3, implies 
$\tilde{C}^2 \ge 3 p_a(\tilde{C}) + K_{\tilde{X}}^2 + \chi (\mathcal{O}_{\tilde{X}}) - 3$, hence we have 
$K_{\tilde{X}} . \tilde{C} \le -(p_a(\tilde{C}) -1) - K_{\tilde{X}}^2 - \chi (\mathcal{O}_{\tilde{X}}) = 
-(g-1) - K_{\tilde{X}}^2 - \chi (\mathcal{O}_{\tilde{X}})$, exactly as in the integral case.
\end{remark}

\begin{remark}
Under the same assumptions, the argument above yields the following stronger inequality: 
$$
K_X.C \le -(g-1) - K_X^2 - \chi (\mathcal{O}_X) - \sum_{p \in C}(\mathrm{mult}_p(C)-1).
$$
\end{remark}

\begin{remark}
The argument above does not apply to the extension to singular curves of item (a) of Theorem \ref{Hao}  
because the analogue of item (iii) does not work. 
Indeed, if $k_{\tilde{C}} \le 4(g - 1) + 3 c_2(\tilde{X}) - K_{\tilde{X}}^2$ then 
$k_C = k_{\tilde{C}} - m \le 4(g - 1) + 3(c_2(X)+1) - (K_X^2 - 1) - m = 4(g - 1) + 3 c_2(X) - K_X^2 + (4-m)$, 
with $4-m > 0$ if $2 \le m \le 3$. 
\end{remark}

\end{document}